\documentclass[12pt]{amsart}

 \usepackage{amsmath}
\usepackage{amssymb,latexsym,cite}
\usepackage{amscd}
\usepackage{comment}
 \usepackage{xspace}
 \usepackage{verbatim}
 \usepackage{todonotes}

\newtheorem{theorem}{Theorem}[section]
\newtheorem{proposition}[theorem]{Proposition}
\newtheorem{corollary}[theorem]{Corollary}
\newtheorem{lemma}[theorem]{Lemma}
\theoremstyle{definition}

\newtheorem{definition}[theorem]{Definition}
\numberwithin{equation}{section}

\newcommand{\R}{{\mathbb R}}

\newcommand{\BSA}{\begin{subarray}}
\newcommand{\ESA}{\end{subarray}}
\newcommand{\BAL}{\begin{aligned}}
\newcommand{\EAL}{\end{aligned}}

\newcommand{\note}[1]{\noindent\textit{#1.}\hspace{2mm}}
\newcommand{\Remark}{\note{Remark}}
\newcommand{\forevery}{\quad \forall}

\newcommand{\rec}[1]{\frac{1}{#1}}

\newcommand{\dist}{\mathrm{dist}\,}

\newcommand{\prt}{\partial}
\newcommand{\sms}{\setminus}
\newcommand{\ems}{\emptyset}
\newcommand{\ti}{\times}

\newcommand{\tl}{\tilde}

\newcommand{\sbs}{\subset}

\newcommand{\nin}{\not\in}

\newcommand{\sth}{such that\xspace}

\newcommand{\bdw}{\partial\Gw}

\newcommand{\qtxt}[1]{\quad\textrm{#1}}
\newcommand{\ntxt}[1]{\noindent\textit{#1}.}

\def\ga{\alpha}     \def\gb{\beta}       \def\gg{\gamma}
       \def\gd{\delta}      \def\ge{\epsilon}
\def\gth{\theta}                         
       \def\vgf{\varphi}    
      \def\gk{\kappa}      \def\gl{\lambda}

      \def\gw{\omega}
\def\gx{\xi}                \def\gz{\zeta}
\def\Gg{\Gamma}     \def\Gd{\Delta}      

\def\Gl{\Lambda}          \def\Gp{\Pi}
\def\Gw{\Omega}              

   \def\BBR {\mathbb R}

\def\gV{{\gg V}}

\begin{document}

\title[Estimates]{Estimates of Green and Martin kernels for Schr\"odinger operators with singular potential in Lipschitz domains}
\author{Moshe Marcus}

\date{\today}

\maketitle

\begin{abstract}
	Consider operators of the form $L^{\gg V}:=\Gd +\gg V$ in a bounded Lipschitz domain $\Gw\sbs \mathbb{R}^N$. Assume that $V\in C^1(\Gw)$ satisfies $|V(x)| \leq \bar a \dist(x,\bdw)^{-2}$ for every $x\in \Gw$ and $ \gg$ is a number in a range $(\gg_-,\gg_+)$ described in the introduction.The model case is $V(x)= \dist(x,F)^{-2}$ where $F$ is a closed subset of $\bdw$ and $\gg< c_H(V)=$ Hardy constant for $V$. We provide sharp two sided estimates of the Green and Martin kernel for $L^\gV$ in $\Gw$. In addition we establish a pointwise version of the 3G inequality.
\end{abstract}

\section{Introduction}

Let $\Gw$ be a bounded Lip domain in $\R^N$, $N\geq3$. We study the operator
$$L^{\gg V}:=\Gd +\gg V$$
where $V\in C^1(\Gw)$ and $\gg$ is a constant. We assume that the
potential $V$ satisfies the conditions:
\begin{equation}\label{Vcon1}
 \exists \bar a>0\, :\quad |V(x)| \leq \bar a \gd(x)^{-2} \forevery x\in \Gw
\end{equation}
$$ \gd(x)=\gd_\Gw(x):=\dist(x,\bdw)$$
and
\begin{equation}\label{Vcon2}
 \gg_-<\gg<\gg_+,
\end{equation}
where
\begin{equation}\label{gg+-}\BAL
\gg_+&=\sup  \{\gg: \exists u_\gg>0\;\text{\sth}\; L^{\gg V}u_\gg=0\},\\
\gg_-&=\inf\,  \{\gg: \exists u_\gg>0\;\text{\sth}\; L^{\gg V}u_\gg=0\}.
\EAL\end{equation}

By a theorem of Allegretto and Piepenbrink \cite{BS} or \cite[Theorem 2.3]{YPrev}, \eqref{gg+-} is equivalent to,
\begin{equation}\label{ggpm}\BAL
\gg_+&=\sup\{\gg: \int_\Gw |\nabla\phi|^2\,dx\geq \gg\int_\Gw \phi^2 V\, dx  \forevery \phi\in H^1_0(\Gw)\},\\
\gg_-&=\inf\,\{\gg: \int_\Gw |\nabla\phi|^2\,dx\geq \gg\int_\Gw \phi^2 V\, dx  \forevery \phi\in H^1_0(\Gw)\}.
\EAL\end{equation}

If $V$ is positive $\gg_+$ is \emph{the Hardy constant relative to $V$} in $\Gw$, denoted by $c_H(V)$. Condition \eqref{Vcon1} and Hardy's inequality imply that $\gg_+>0$ and $\gg_-<0$. Clearly,
if $V>0$ then $\gg_-=-\infty$ and if $V<0$ then $\gg_+=\infty$. Finally, if $\gg\in (\gg_-,\gg_+)$ then there exists a Green function for  $L^{\gg V}$ in $\Gw$, denoted by, $G^{\gg V}_\Gw$. The subscript will be dropped, except when several domains are considered.

Assumptions \eqref{Vcon1} -- \eqref{Vcon2} imply that $-L^{\gV}$ is positive and its first eigenvalue
$\gl_\gV$
is positive. The corresponding normalized eigenfunction is denoted by $\vgf_\gV$. (The normalization is
$\vgf_\gV(x_0)=1$ where $x_0$ is a fixed reference point in $\Gw$.)

 The following result is due to Pinchover \cite{YP-per}. It is proved by adapting an argument from \cite[Theorem 3]{YP89}.

 \begin{lemma}\label{YP1} Assuming \eqref{Vcon1} -- \eqref{Vcon2}, there exists $\ge>0$ \sth the operator $-(L^{\gg V} +\ge\gd(x)^{-2})$ has a positive supersolution.
 \end{lemma}
For the convenience of the reader, a  proof is provided in the next section.

This fact implies that $L^{\gg V}$, $\gg\in (\gg_-,\gg_+)$, is \emph{weakly coercive} in the sense of \cite{An87}. Therefore one may apply to it potential theoretic results of
 Ancona \cite{An87} and \cite{An97}. In particular one may apply to this operator the Boundary Harnack Principle \cite{An87}, that plays a crucial role in the present work.
\medskip

\note{Notation} Let $f,g$ be non-negative functions in a domain $D$. The notation $f\sim g$ in  $D$ means that there exist two positive constants $c_1,c_2$ -- called \emph{similarity constants} --
 \sth
$$c_1f\leq g\leq c_2f \qtxt{in }D.$$

 The notation $f\lesssim g$ in  $D$ means that there exist a positive constant $c$ \sth
$$f\leq cg \qtxt{in }D.$$ 


 \begin{lemma}\label{AA1} Assume \eqref{Vcon1} -- \eqref{Vcon2}. For any $x_0\in \Gw$ and $\ge>0$,
\begin{equation}\label{G/vgf}
  G^\gV(\cdot, x_0) \sim \vgf_\gV(\cdot) \qtxt{in }\{x\in\Gw:\, |x-x_0|>\ge\}.
\end{equation}
Of course, the similarity constants depend on $x_0$ and $\ge$.
\end{lemma}

\Remark The fact that \eqref{G/vgf} holds for every potential satifying \eqref{Vcon1} and \eqref{Vcon2} was pointed out to me by Alano Ancona.

For the convenience of the reader, a proof is provided in the next section.

  It is well known  that, in every compact set $F\sbs \Gw$,
\begin{equation}\label{GinF}
  c_1(F)|x-y|^{2-N}\leq G^\gV(x,y)\leq c_2(F)|x-y|^{2-N}. 
\end{equation}
 Sharp,  two sided estimates    of the Green kernel of the Laplacian, \emph{up to the boundary,} in Lipschitz domains, were obtained by Bogdan \cite{Bog}. In smooth domains such estimates have been obtained in \cite{FMT}, when $V=\gd^{-2}$  (= the classical Hardy potential) and  $\gg\in (0, c_H(V))$. These estimates can be extended to a large class of potentials using results on comparison of Green functions for related operators.
 In the case of  \emph{small perturbations} of the potential results of this type  were obtained by Murata \cite{Mur86}, \cite{Mur97} and Pinchover \cite{YP89},  \cite{YP90}. See \cite{YPrev} for a survey of these and related papers. These results are obtained without any explicit assumptions on the domain, which may also be unbounded.  However the assumptions on the operators -- including the  existence of a positive minimal Green function -- and the definition of a `small' perturbation reflect implicitly on the domain.

The results of Ancona \cite{An87} imply the existence of the Green function for a large class of potentials in bounded Lipschitz domains and even more general cases (e.g. John domains).

Combining the results of \cite{YP90} with those of \cite{An87} one obtains for instance the following:

  Assume that $\Gw$ is a bounded Lipschitz domains and that $V\in C^\ga(\Gw)$, $\ga\in (0,1]$, and $\gg$ satisfy conditions \eqref{Vcon1} and \eqref{Vcon2}.  Let $W:=\gg V+V_0$ where $V_0\in C^\ga(\Gw)$ and
$$|V_0|\leq c\gd^{\ge-2} \qtxt{in}\;\Gw$$
for some positive numbers $c,\ge$. Then
 $$G^\gV_\Gw \sim G^W_\Gw.$$
  In particular, letting $V=0$, we conclude that, for $V_0$ as above, the Bogdan estimates hold for the Green kernel of the operator $-\Gd + V_0$ in bounded Lipshitz domains.
\smallskip

 This is also  a consequence of the results of  Ancona \cite{An97} in which the author established the equivalence of the Green functions for a pair of operators $L^{V_i}$, $i=1,2$  under very general conditions on $V_1-V_2$.
 In fact the results of \cite{An97} apply to more general Schr\"odinger operators, where $\Gd$ is replaced by a  linear second order  elliptic operators whose coefficients may be singular on $\bdw$.  In this case the conditions are imposed on the weighted difference of the Schr\"odinger operators.

Sharp estimates have also been obtained for Hardy potentials in conical domains, possibly unbounded,   e.g. \cite{DPP}, \cite{Hirata}.

In the case of smooth domains and potentials with \emph{singularities in $\Gw$}, two sided estimates of the Green function have been obtained under very general conditions,  see \cite{F+VI}, \cite{FN+VI}, \cite{VI} and references therein.
These estimates are sharp with respect to $\ln G$.

\medskip

In the present paper we derive sharp, up to the boundary,  two-sided estimates of the Green kernel of $L^{\gV}$ in bounded Lipschitz domains.  
Following are the main results.

\begin{theorem}\label{G-estA} Assume \eqref{Vcon1} -- \eqref{Vcon2} and $N\geq 3$.

 Then, for every $b>0$ there exists
a constant $C(b)$, depending also on $N, r_0, \kappa, \bar a$, \sth: if $x,y\in \Gw$ and
\begin{equation}\label{d>>r}
 |x-y|\leq \rec{b}\min(\gd(x),\gd(y))
\end{equation}
then
\begin{equation}\label{G-est.1}\BAL
 \rec{C(b)}|x-y|^{2-N}\leq G^\gV_\Gw(x,y)\leq C(b)|x-y|^{2-N}.
\EAL\end{equation}
\end{theorem}
\medskip

In the next theorems, $C$ stands for a constant depending only on $r_0,\kappa, \bar a$ and $N$.

\begin{theorem}\label{G-estB} Assume \eqref{Vcon1} -- \eqref{Vcon2} and $N\geq 3$.

 If $x,y\in \Gw$ and
\begin{equation}\label{d<r0}
\max(\gd(x),\gd(y))\leq r_0/10\kappa
\end{equation}
\begin{equation}\label{d<<r}
 \min(\gd(x),\gd(y))\leq \frac{|x-y|}{16(1+\kappa)^2}
\end{equation}
then
\begin{equation}\label{G-est.2}\BAL
\rec{C}|x-y|^{2-N}\frac{\vgf_\gV(x)\vgf_\gV(y)}{\vgf_\gV( x_y)\vgf_\gV( y_x)}&\leq G^\gV_\Gw(x,y)\\
&\leq
C|x-y|^{2-N}\frac{\vgf_\gV(x)\vgf_\gV(y)}{\vgf_\gV( x_y)\vgf_\gV( y_x)}.
\EAL\end{equation}
The points $x_y$, $y_x$ depend on the \emph{pair} $(x,y)$. If
$$\hat r(x,y):=|x-y|\vee \gd(x)\vee \gd(y)\leq r_0/10\kappa$$
they can be chosen arbitrarily in the set
\begin{equation}\label{Axy}
 A(x,y):=\{z\in \Gw: \rec{2}\hat r(x,y)\leq \gd(z)\leq 2\hat r(x,y)\} \cap B_{4\hat r(x,y)}(\frac{x+y}{2})\}.
\end{equation}
Otherwise set $ x_y= y_x= x_0$ 
where $x_0$ is a fixed reference point.
\end{theorem}

\Remark  There exists a constant $C$ \sth for any two points $x,y\in\Gw$ and any $P,Q\in A(x,y)$,
$$\rec{C}\vgf_\gV(P)\leq \vgf_\gV(Q)\leq C\vgf_\gV(P).$$
This is a consequence of the strong Harnack inequality (see Lemma \ref{b<b_0} below)
and the fact that, under condition \eqref{d<<r}, $\hat r(x,y)\sim |x-y|$.  

The same observation is valid if $A(x,y)$ is replaced by
\begin{equation}\label{Abxy}
A_b(x,y):=\{z\in \Gw: \rec{b}\hat r(x,y)\leq \gd(z)\leq b\hat r(x,y)\} \cap B_{4\hat r(x,y)}(\frac{x+y}{2}) 
\end{equation}
where $b$ is a number in $(1,r_0/10b\kappa)$ and $C$ is a constant depending on $b$. Consequently Theorem \ref{G-estB} remains valid if $A(x,y)$ is replaced by $A_b(x,y)$ and $C$ by $C_b$, i.e. a constant depending on $b$.
\medskip

Let $K^\gV_\Gw$ denote the Martin kernel of $L^\gV$ in $\Gw$. As a consequence of the previous result we obtain

\begin{theorem}\label{K-est} Assume \eqref{Vcon1} -- \eqref{Vcon2} and $N\geq 3$.

If $x\in \Gw$, $y\in \bdw$ and $|x-y|< \frac{r_0}{10\kappa}$ then
\begin{equation}\label{K-est.1}\BAL
\rec{C} \frac{\vgf_\gV(x)}{\vgf_\gV(x_y)^2}|x-y|^{2-N}\leq K_\Gw^\gV(x,y)\leq
C \frac{\vgf_\gV(x)}{\vgf_\gV(x_y)^2}|x-y|^{2-N},
\EAL\end{equation}
where $x_y$ is an arbitrary point in $A(x,y)$.
\end{theorem}
\smallskip

\begin{definition}
Let $\gz\in \bdw$. A unit vector $\nu$ in $\mathbb{R}^N$ is an \emph{inner pseudo normal} at $\gz$ if 
$$a_\nu(\gz):=\limsup_{x\in \bdw;\,x\to\gz} \frac{\langle x-\gz,\nu\rangle}{|x-\gz|}<1.$$
Let $\gl\in (0,1)$. The vector $\nu$ is a \emph{$\gl$ -- inner normal} at $\gz$ if $a_\nu(\gz)<\gl$.

\end{definition}

Another consequence of the previous estimates is the following version of the 3G inequality. 
\smallskip

\begin{theorem}\label{3G} Assume \eqref{Vcon1} -- \eqref{Vcon2} and $N\geq 3$. In addition assume that there exist numbers  $\Gl\in (0,1)$ and $b_1>1$ \sth, if $x,y\in \Gw$ satisfy \eqref{d<r0} and $x,y$ lie on a $\gl$ -- inner normal at $\gz$ for a point $\gz\in \bdw$ and some $\gl<\Gl$   then  
\begin{equation}\label{q-mon}\BAL
 b_1\gd(y)\leq \gd(x)
\Longrightarrow \vgf_\gV(y)\leq  C(b_1)\vgf_\gV(x).
\EAL\end{equation}

Under these assumptions, if $x,y,z$ are three distinct points in $\Gw$ then,
\begin{equation}\label{3G.1}
 \frac{G^\gV(x,y)G^\gV(y,z)}{G^\gV(x,z)}\leq C'(b_1)\big( |x-y|^{2-N} + |y-z|^{2-N} \big).
\end{equation}
\end{theorem}

\Remark The 3G inequality, in various forms, has been studied in numerous papers.  In  \cite{An97} it was established with respect to a  larger class of potentials -- in particular, without assuming \eqref{q-mon} --  but with certain restrictions on the configuration of the three points. 

\section{Notations and preliminaries}

We start with the proof of two auxiliary lemmas stated in the introduction.
\smallskip

\note{Proof of Lemma \ref{YP1}} Let $V_1, V_2\in C^2(\Gw)$ and assume that $V_i$ satisfies \eqref{Vcon1} and that $\gg_-(V_i)<1<\gg_+(V_i)$, $i=1,2$. The latter assumption implies the existence of a positive eigenfunction, say $u_i$,  of $-L^{V_i}$, with positive eigenvalue. Let $a\in (0,1)$ and put,
$$u :=u_1^a u_2^{1-a}, \quad W:=aV_1+(1-a)V_2.$$
A straight forward computation yields:

\begin{equation}\label{comp1}
 -L^W(u)=-\Gd u- Wu = f_1+f_2
\end{equation}
where
$$\BAL
f_1 &= -\frac{au}{u_1}L^{V_1}u_1 - \frac{(1-a)u}{u_2}L^{V_2}u_2 \geq 0\\
f_2 &= a(1-a)u \left[\frac{\nabla u_1}{u_1}-  \frac{\nabla u_2}{u_2}\right]^2\geq 0.
\EAL$$
Thus \emph{$u$ is a positive $L^{W}$~superharmonic function.}
\smallskip

Now we apply this result to the following case: \\[3mm]
$V_1=\gg V$ where $V$ and $\gg$ satisfy \eqref{Vcon1} and \eqref{Vcon2},\\[2mm]
$V_2= \dfrac{c_H}{2} \gd^{-2} $ where $c_H$ is the classical Hardy constant in $\Gw$.
\smallskip

It follows that for every $a\in (0,1)$, the function $v=\vgf_\gV^a\vgf_0^{1-a}$ $L^W$~superharmonic in $\Gw$, i.e.,
$$-(\Gd + a\gg V)v \geq \ge \gd^{-2},  \quad \ge:=\frac{(1-a)c_H}{2}.$$
As this result is valid for any $\gg\in (\gg_-,\gg_+)$ and any $a\in (0,1)$, Lemma \ref{YP1} holds.

\qed

\note{Proof of Lemma \ref{AA1}} Put $V_1=\gV$ and $V_2=\gV + \gl_{\gV}$. Then $\vgf_{\gV}$ is a minimal positive ground state solution of $L^{V_2}$. Evidently $W:=V_2-V_1$ is  a small perturbation of $V_1$. Therefore \eqref{G/vgf} is a consequence of Theorem 3.1  and Lemma 3.6 of \cite{YP89}.

\qed

Next we introduce some notations that will be used throughout the paper.

Given $r,\rho$ positive  denote
$$ T^0(r,\rho)=\{\xi=(\xi_1,\xi')\in \R\ti\R^{N-1}: |\xi'|<r, \; |\xi_1|< \rho\}.$$
If $\xi=\xi^P$ is centered at $P$ we denote by
$T^P_{\xi^P}(r,\rho)$ the cylinder  $T^0(r,\rho)$ in this set of coordinates. However, as a rule we shall drop the subscript ${\xi^P}$.

Since $\Gw$ is a bounded Lipschitz domain there exists $\kappa\geq 1$ and $r_0>0$ \sth, for every
$P\in \bdw$, there exists an Euclidean set of coordinates $\xi=\xi^P$, centered at $P$, and a
$\kappa$-Lipschitz function $f^P:\R^{N-1}\mapsto \R$ \sth $f^P(0)=0$ and
\begin{equation}\label{TP}
  T^P(r_0,10\gk r_0)\cap\Gw=\{(\xi_1,\xi'): |\xi'|<r_0, \; f^P(\xi')<\xi_1<10\gk r_0\}
\end{equation}
in the set of coordinates $\xi^P$. Any set of coordinates centered at $P$ \sth \eqref{TP} holds is called an \textit{admissible set of coordinates at $P$} and $ T^P(r,\rho)$, $r\in (0,r_0)$ and $\rho \in (0, 10\kappa r_0]$ is called a \textit{standard cylinder at $P$.}  The couple $(\gk, r_0)$ is called the
\textit{Lipschitz characteristic} of $\Gw$. It is not unique, but will be kept fixed throughout the paper.

For $r\leq r_0$, $\rho\leq 10\kappa r_0$ we denote
\begin{equation}\label{gwP}
  \gw^P(r,\rho):= T^P(r,\rho)\cap\Gw
\end{equation}
where $T^P(r,\rho)$ is  a
\textit{standard cylinder} at $P$.  If $\xi\in \gw^P(r_0, 10\kappa r_0)$, $\xi_1>0$ and $|\xi'|/\xi_1<\kappa/2$ we say
that the unit vector in the direction
$\stackrel{\longrightarrow}{P\xi}$ is an \textit{approximate normal at} $P$. This vector is denoted by $\mathbf{n}(P,\xi)$.

\medskip

The boundary Harnack principle (briefly BHP) due to \cite{An87} plays a crucial role in
the paper. For easy reference, we state it below.

\begin{theorem}\label{BHP} Let $P\in \bdw$ and let $T^P(r,\rho)$ be a standard cylinder at $P$.
      There exists a constant $c$ depending only on $N,\, \bar a$ and ${ \frac {   \rho } { r}}$ such
      that
      whenever $u$ is a positive $L^\gV$ harmonic function  in $\omega^P(r,\rho)$  that vanishes
      continuously
      on $\bdw \cap T^P(r, \rho )$
 then
\begin{equation}\label{u/v+}
\hspace{2truemm}  c^{-1}r^{N-2}\, \mathbb G ^{V}_\Gw(x,A')\leq \frac{u(x)}{u(A)}\leq c\, r^{N-2}
\,\mathbb G
^{V}_\Gw(x,A'), \quad   \forall x\in \Omega  \cap \overline T^P({ \frac {  r} {2 }}; {
\frac {\rho } {2 }})
\end{equation}
where $A=(\rho/2)(1,0,...,0)$, $ A'=(2\rho/3)(1,0...,0)$ in the corresponding set of local coordinates
$\xi^P$.

In particular,  for any pair $u,v $ of  positive $L^\gV$ harmonic functions in  $ \omega^P(r,\rho) $ that
vanish on
$\bdw \cap T^P(r, \rho )$:
\begin{equation}\label{u/v}
 u(x)/v(x)\le Cu(A)/v(A), \quad \quad  \forall x\in \Gw\cap\overline T^P(r/2, \rho /2))
\end{equation}
where $C=c^2$.
 \end{theorem}

\Remark (i) Inequality \eqref{u/v+} implies that (in the notation of the theorem)
\begin{equation}\label{eAA'}
 \rec{c}r^{2-N}\leq G^\gV_\Gw(A,A')\leq cr^{2-N}.
\end{equation}
(ii) Inequality \eqref{u/v+} remains valid for any $A, A'$ \sth $A=(a_1\rho,0...,0)$, $A'=
(a_2\rho,0...,0)$ and $0<a_1<a_2<1$. In this case, the constant $c$  depends also on $a_1/ a_2$.

\section{Proof of Theorem \ref{G-estA}}

The proof is based on several lemmas in which we assume, without further mention,  that conditions
\eqref{Vcon1} -- \eqref{Vcon2} are satisfied.
\medskip

\noindent\textit{Notation} (i) Put
\begin{equation}\label{Ggy}
 \Gg_y(x):=a_N|x-y|^{2-N}
\end{equation}
where $a_N$ is the  constant \sth $-\Gd
\Gg_y=\gd_y$.

(ii) Denote by
$G^\gV_y$ the function $x\mapsto G^\gV(x,y)$.

(iii) For every $b>1$ and $z\in \Gw$, put
\begin{equation}\label{Byb}\BAL
 B^z_b &=\{x\in \Gw:\, |x-z|<\frac{\gd(z)}{b}\}, \\   S^z_b &= \{x\in \Gw:\, |x-z|=\frac{\gd(z)}{b}\}.
\EAL\end{equation}

\begin{lemma}\label{AA'1}
Let $z\in \Gw$, $\gd(z)<r_0/4$. For every $b>1$ there exists a constant $c^*=c^*(b)>0$, independent of $z$,
\sth
\begin{equation}\label{eAA'1}
\rec{c^*}\gd(z)^{2-N}\leq G^\gV_\Gw(x,z)\leq c^*\gd(z)^{2-N} \forevery x\in S^z_b
\end{equation}
\end{lemma}

\proof First we prove,

\noindent\textsc{Assertion 1.} \  \textit{For every $z$ as above, there exists a constant $c_1$ depending on $b$ but independent of $z$  and a point $\gz\in S^z_b$ \sth}
\begin{equation}\label{eAA'3}
 \rec{c_1}\gd(z)^{2-N}\leq G^\gV_\Gw(\gz,z)\leq c_1\gd(z)^{2-N}.
\end{equation}

Let $Q\in \bdw$ be a point \sth $|Q-z|=\gd(z)$. Let $T^Q(r_0,10\kappa r_0)$ be a standard cylinder at $Q$
associated with a local set of coordinates $\xi^Q$.

Let $P$ be the point on $\bdw$ \sth $(\xi^Q)'(P)=(\xi^Q)'(z)$. Then
$\xi^{P}:=\xi^Q-\xi^Q(P)$ is a local set of coordinates at $P$ and
$T^{P}(r_0/2, 5r_0\kappa)$ is a standard cylinder at $P$ relative to $\xi^P$. (Recall that
$\gd(z)<r_0/4$.)

Let $\gz$ be the point of intersection of the segment $[P,z]$ with the sphere $S^z_b$. We apply Theorem
BHP in $T^{P}(2\gd(z), 20\kappa\gd(z))$ when $A'=z$ and $A=\gz$. This is possible because
$$\BAL \gd(z)&\leq|Pz|<|PQ| +|Qz|\leq (\kappa +1) \gd(z),\\
\frac{b-1}{b}\gd(z)&\leq|P\gz|=|Pz|-\gd(z)/b\leq(\kappa +1 -\frac{1}{b})\gd(z)
\EAL$$
and consequently,
$$1-\rec{b}\leq \frac {|P\gz|}{|Pz|}\leq 1 -\rec{b(\kappa+1)}.$$
Thus \eqref{eAA'3} - with a constant $c_1$ depending on $b$ but independent of $z$ - is a consequence of
\eqref{eAA'} and the remark following it.

\medskip

  The sphere $S^z_b$ can be covered by $c'(N)$ balls of radius
$r'=\gd(z)/4b$
centered on the sphere. If $x\in S^z_b$ then $\gd(x)\geq \frac{b-1}{b}\gd(z)\geq 4(b-1)r'$. Therefore
by the classical Harnack inequality, there exists a constant $C'(N)$ (independent of $z$) \sth
\begin{equation}\label{sup<inf}
 \sup_{x\in S^z_b} (G^\gV_\Gw)(x,z)  \leq C'(N)\inf_{x\in S^z_b} (G^\gV_\Gw)(x,z).
\end{equation}
This inequality and Assertion 1 imply \eqref{eAA'1}.

\qed
\medskip

\begin{lemma}\label{b<b_0} Let $F\in C(\Gw)$ be a positive function satisfying the strong Harnack inequality. Let $b,b_0$ be two numbers \sth $0<b<b_0$, let $x,y\in \Gw$ and put $r=|x-y|$. Suppose that
\begin{equation}\label{eAA'4}\BAL
r<\frac{r_0}{10\kappa b_0},\quad br&\leq \min(\gd(x),\gd(y))\\
&\leq \max(\gd(x),\gd(y))\leq (b_0+1)r,
\EAL\end{equation}
where $(\kappa,r_0)$ is the Lipschitz characteristic of $\Gw$ (see Section 2).
Then there exists a constant  $c^*$, independent of $x,y$ (but depending on $N, \kappa, r_0, b, b_0$ and the Harnack constants for $F$)  \sth
\begin{equation}\label{eAA'0}
\rec{c^*}F(x)\leq  F(y)\leq c^*F(x).
\end{equation}
\end{lemma}

\proof  Let $X\in \bdw$ be a point \sth $|x-X|=\gd(x)$. Let $\gx^X$ be an admissible set of local coordinates at
$X$ associated with the cylinder $T^X(r_0,10\kappa r_0)$ (see \eqref{TP}). Put
$$r_1=|\xi^X_1(x-y)|, \quad r'= |(\xi^X)'(x-y)|.$$
Let $Y\in\bdw$ be the point  \sth $(\xi^X)'(Y)=(\xi^X)'(y)$ and let
$\xi^Y$ be the set of coordinates centered at $Y$ given  by
$$\xi^Y=\xi^X - \xi^X(Y).$$

Denote
$$\tl\gd(z):=\xi_1^X(z)- f^X((\xi^X)'(z)) \forevery z\in  T^X(r_0,10\gk r_0)\cap\Gw.$$
Note that $\tl\gd(z)$ is simply the distance between $z$ and $\bdw$ measured along the line through $z$ parallel to the $\gx_1$ axis.
The definition of standard cylinder (see \eqref{TP}) implies that          $\tl\gd(z)>0$ for every $z$ as above. In addition,
\begin{equation}\label{tld}
\frac{\tl\gd(z)}{\sqrt{1+\kappa^2}}\leq \gd(z)\leq \tl\gd(z).
\end{equation}
The right inequality is trivial and the left follows from the Lipschitz property of $\bdw$ and the previous remarks on  $\tl\gd$.

Put $d=\tl\gd(x)$ and
\begin{equation}\label{def_J}\BAL
J&=(\bdw\cap \overline{T}^X(r_0/2,10\kappa r_0))+(d,0,\dots,0)\\
&=\{z:\, \xi^X_1(z)=f^X((\xi^X)'(z))+d,\; |(\xi^X)'(z)|\leq r_0/2\}.
\EAL\end{equation}
Note that $x\in J$ (but $y$ need not be in $J$) and
\begin{equation}\label{tldJ}
\tl \gd(z) =d \forevery z\in J.
\end{equation}

Denote by $y^*$ the $\xi^X_1$-projection of $y$ on $J$:
$$\xi^X_1(y^*)= f^X((\xi^X)'(y))+d,\quad (\xi^X)'(y^*)=(\xi^X)'(y).$$
If $[x,y]$ is parallel to the $\xi_1$ axis, it is easy to see that \eqref{eAA'0} holds. Therefore we may assume that $x\neq y^*$. 

Let $\Gp$ denote the plane containing $x,y$ that is parallel to the $\xi^X_1$ axis. Then $y^*\in \Gp$ and we denote by $I_J(x,y^*)$ the closed section of the curve $\Gp\cap J$ with end points $x, y^*$. Let $I(x,y)$ be the curve connectiong $x,y$ given by,
\begin{equation}\label{Ixy1}
I(x,y)= I_J(x,y^*)\cup[y^*,y].
\end{equation}
If $y\neq y^*$ and $\gth$ is the angle between $I_J(x,y^*)$ and $[y^*,y]$ then $|\cot\gth|\leq \rec{k}$. Therefore $I(x,y)$ is a  Lipschitz curve.

By definition, $\gd(x)\leq d=\tl\gd(x)$. Therefore, by \eqref{eAA'4},  and \eqref{tld},
\begin{equation}\label{d-est}
br  \leq d \leq \gd(x)\sqrt{1+\kappa^2}\leq (b_0+1)r\sqrt{1+\kappa^2}.
\end{equation}
By \eqref{tld} and \eqref{tldJ},
\begin{equation}\label{gdz}
\frac{d}{\sqrt{1+\kappa^2}}\leq \gd(z)\leq d \forevery z\in J.
\end{equation}
Hence,
\begin{equation}\label{gdJ}
\gd(J):=\min_J\gd(z)\geq \frac{br}{\sqrt{1+\kappa^2}}.
\end{equation}

For every $z\in [y^*,y]$, $(\xi^X)'(z)= (\xi^X)'(Y)$. Therefore
$\tl \gd(z)$ lies between $\tl \gd(y)$ and $\tl \gd(y^*)$ 
for every $z\in [y^*,y]$. Since $\tl\gd(y)\geq \gd(y)$ and, by \eqref{tldJ},
$\tl \gd(y^*)=d$ it follows that
$$\min_{[y^*,y]}\tl\gd(z)\geq \min(d,\gd(y)).$$
Hence, by \eqref{eAA'4}, \eqref{tld} and \eqref{d-est},
\begin{equation}\label{y*y}
\min_{[y^*,y]}\gd(z)\geq \frac{br}{\sqrt{1+\kappa^2}}.
\end{equation}

The curve $I_J(x,y^*)$ is given by, 
\begin{equation}\label{xi(t)}\BAL
\{\gx^X(t):\; &\gx_1^X(t)=f^X((\gx^X)'(t))+d,\\ &(\gx^X)'(t)=(1-t)(\xi^X)'(x)+t (\xi^X)'(y), \;
 t\in (0,1)\}
\EAL\end{equation}
Therefore
\begin{equation}\label{|xy|J}
|I_J(x,y^*)|\leq \int_{0}^1|\nabla f^X\, \frac{d(\gx^X)'    }{dt}|dt\leq \kappa |(\xi^X)'(x-y)|=\kappa r'.
\end{equation}
Furthermore, by \eqref{eAA'4}, \eqref{tld} and \eqref{d-est}
\begin{equation}\label{y-y*}\BAL
|y-y^*|&=|\tl\gd(y)-\tl\gd(y^*)|=|\tl\gd(y)-d|\leq \max(\tl\gd(y),d)\\
&\leq \max(\gd(y)\sqrt{1+\kappa^2},d) \leq (b_0+1)\sqrt{1+\kappa^2}\,r.
\EAL\end{equation}
Thus the curve $I(x,y)$ has total length no larger than
$C_{b_0} r$ where $C_{b_0}=\kappa + (b_0+1)\sqrt{1+\kappa^2}$

Let $P_0, \cdots, P_m$ be distinct points on $I(x,y)$,  $P_0=x$, $P_m=y$ and let $D_i$ be the
open ball of radius $s:=br/4\sqrt{1+\kappa^2}$ centered at $P_i$, $i=0,\cdots,m$. We assume that the points
$P_i$ are so distributed that
$$D_i\cap D_{i+1}\neq \ems, \quad D_i\cap D_j=\ems\;\mathrm{\ if\ }\; |i-j|>1.$$
 By \eqref{gdJ} and \eqref{y*y}, $\gd(P_i)\geq 2s  $. Since the total
length of  $I_J(x,y^*)\cup [y^*,y]$ is not larger than $C_{b_0}r$, the number of points $m+1$ needed in
order to achieve such a configuration depends only on $b$, $b_0$ and $\kappa$. Therefore, as $F$ satisfies the strong Harnack inequality, \eqref{eAA'0} follows.

\qed

\begin{lemma}\label{b<1} Let $N\geq 3$. Assume that there exists $b_0>1$ \sth the statement of Theorem \ref{G-estA} is valid when $b > b_0$.
Then it is also valid when $b\in (0,b_0]$:

If $0<b\leq b_0$, $x,y\in \Gw$ and
\begin{equation}\label{AA'4}
r:=|x-y|<r_0/10\kappa\,b_0, \quad br\leq  \min(\gd(x),\gd(y))
\end{equation}
then
\begin{equation}\label{AA'5}
\rec{c^*}r^{2-N}\leq G^\gV_\Gw(x,y)\leq c^*r^{2-N},
\end{equation}
where $c^*(b, b_0)$ is a constant  independent of $x,y$.
\end{lemma}

\proof If $\gd(y)>(b_0+1)r$ then $\gd(x)>b_0r$  so that \eqref{d>>r} holds for some $b>b_0$ and therefore
\eqref{G-est.1} holds by assumption. The statement is symmetric in $(x,y)$ so that we may assume:
\begin{equation}\label{max<r}
 \max(\gd(x),\gd(y))\leq (b_0+1)r.
\end{equation}

As  shown in the proof of Lemma \ref{b<b_0}, \eqref{AA'4} and \eqref{max<r} imply that the points $x,y$ can be joined by a Lipschitz curve $I(x,y):= I_J(x,y^*)\cup [y^*,y]$ (notation as in that lemma) \sth:
\begin{equation}\label{Ixy}\BAL
\mathrm{length}\,I(x,y)& \leq C_{b_0}r,\qquad  C_{b_0}=\kappa+(b_0+1)\sqrt{1+\kappa^2}\\
\min_{z\in I(x,y)}\gd(z)& \geq br \sqrt{1+\kappa^2}.
\EAL\end{equation}

Let $t:=br/4b_0\sqrt{1+\kappa^2}$. Then
\begin{equation}\label{2b_0t}
\gd(z)> 2b_0t, \forevery z \in I(x,y).
\end{equation}
Since $r>t$,  $x \not\in B_t(y)$.
 Let $\eta$ be the closest point to $x$ among all points $z\in I(x,y)$ \sth 
 $|z-y|=2t$.  By the assumption of the lemma and \eqref{2b_0t},
\begin{equation}\label{eta1}
\rec{c} t^{2-N}\leq G^\gV_y(\eta)\leq c t^{2-N}.
\end{equation}
Let $I(x,\eta)$ denote the part of $I(x,y)$ connecting $x$ and $\eta$. 
By Lemma \ref{b<b_0} applied to the function $F:=G_y^\gV$ in the domain 
$\Gw\sms B_{t}(y)$ we obtain
\begin{equation}\label{eta2}
\rec{c'}G_y^\gV(x)\leq G_y^\gV(\eta)\leq c' G_y^\gV(x).
\end{equation} 
Here $c'$ depends on $b.b_0,\kappa$ and also on the constant associated with the strong Harnack inequality for $G^\gV_y$ in the domain $\Gw\sms B_{t}(y)$. This constant is independent of $y$ and, for balls $B_t(z)$, it is  independent of $z\in I(x,\eta)$  provided that  \eqref{AA'4} and 
\eqref{max<r} hold. This is a simple consequence of the Boundary Harnack principle.
Finally \eqref{eta1} and \eqref{eta2} imply \eqref{AA'5}.

\qed

%
%
%

  The next result is classical. We list it for easy reference.

\begin{lemma}\label{I.1}   For every $y\in \Gw$,
$$\lim_{x\to y}\frac{G^\gV_y}{\Gg_y}=1.$$
\end{lemma}

\begin{lemma}\label{I.2} Assume $N>3$.
 Let $\gg\in (\gg_-, \gg_+)$ and $y\in \Gw$. For every  $x\in \Gw$, denote
\begin{equation}\label{g_y}
g_y(x):=|x-y|^{3-N}\frac{\vgf_\gV(y)}{\gd(y)}
\end{equation}
Then,   there exists a number $b_0>1$  dependent on $\gg V$ and $\bar a$, but not on $y$, \sth, for every
$b>b_0$ there exists a constants $c>0$, dependent on $b$ but not on $y$, \sth
\begin{equation}\label{I.2.1}\BAL
 -L^\gV(\vgf_\gV\Gg_y+cg_y)\geq \vgf_\gV \gd_y  \quad \textrm{in } B^y_b.
\EAL \end{equation}
\end{lemma}

\proof Let $b>1$. By the strong Harnack inequality,
\begin{equation}\label{H-vgf}
 \frac{ \sup_{B^y_b}\vgf_\gV}{\inf_{B^y_b}\vgf_\gV}= c'(b)<\infty.
\end{equation}
A straightforward computation yields,
\begin{equation}\label{Lfg} \BAL
-L^\gV (\vgf\Gg_y)=\, & \vgf(x)\gd_y - 2\nabla \Gg_y(x)\,\cdot\, \nabla\vgf(x)\\
 & - \Gg_y(x)(\Gd\vgf(x) -\gg V\vgf(x))\\
=\,& \vgf(x)\gd_y - 2\nabla \Gg_y(x)\,\cdot\, \nabla\vgf(x) + \gl_1\Gg_y(x) \vgf(x).
\EAL
\end{equation}

By interior elliptic
 estimates (see e.g. \cite[Theorem 6.2]{GT}) and \eqref{H-vgf},

\begin{equation}\label{gradf}
|\nabla\vgf(\xi)|\leq C_1\rec{\gd(\xi)}\sup_{|\xi-x|<\gd(\xi)/b}\vgf(x)\leq
C_2\frac{\vgf(\xi)}{\gd(\xi)} \forevery
\xi\in \Gw.
\end{equation}
The constant $C_2$ is independent of $y$.
By \eqref{H-vgf} and \eqref{gradf},
\begin{equation}\label{h_y}\BAL
|2\nabla \Gg_y(x)\,\cdot\, \nabla\vgf(x)|&\leq C_0|x-y|^{1-N}\frac{\vgf(y)}{\gd(y)}
 =:h_y(x) \forevery x\in B^y_b.
\EAL\end{equation}
The constant $C_0$ is independent of $y$.

By \eqref{Lfg},
\begin{equation}\label{I.2.3}\BAL
-L^\gV (\vgf_\gV\Gg_y)\geq \vgf_\gV\gd_y - h_y, \forevery x\in B^y_b.
\EAL\end{equation}

\medskip

Denote
\begin{equation}\label{Lmu}
 L_\mu:=\Gd+ \frac{\mu}{\gd^2}
\end{equation}
and let $\mu= \gg\bar a$ for $\bar a$ as in \eqref{Vcon1}. Then for any positive  function $f\in L^1_{loc}(\Gw)$ and any $\gg>0$
\begin{equation}\label{LVf}
-L_\mu f\leq -L^\gV f, \quad  -L_{-\mu} f\leq -L^\gV f
\end{equation}
 The second inequality is valid because we assume $|V|\leq \bar a\gd^{-2}$.

If $f_y(x):=|x-y|^{3-N}$, a simple calculation yields
\begin{equation}\label{I.2.4=}
 -L_\mu f_y= \big(N-3 -\mu \frac{|x-y|^2}{\gd(x)^2}\big)|x-y|^{1-N}.
\end{equation}
$$$$
For $x\in B^y_b$, $|x-y|\leq \gd(y)/b$ and $(1-\rec{b})\gd(y)<\gd(x)$ so that
$$|x-y|^2/\gd(x)^2<(b-1)^{-2}.$$
 Therefore,
\begin{equation}\label{I.2.4}
  -L_\mu f_y\geq \ell |x-y|^{1-N},\qquad \ell:=N-3-\frac{\mu}{(b-1)^2}.
\end{equation}

If $N>3$, let $\mu=\gg\bar a$ and let $b_0$ be sufficiently large so that $\ell>0$ for $b\geq b_0$. Note that the choice of $b_0$ does not depend on $y$.

Pick a constant $c$  \sth
 $$c\ell>C_0, \quad \textrm{$C_0$ as in \eqref{h_y}}.$$
Then  by \eqref{g_y}, \eqref{h_y}, \eqref{I.2.3}(a),
\eqref{LVf} and \eqref{I.2.4},
\begin{equation}\label{I.2.5}\BAL
 -L^\gV(\vgf\Gg_y+cg_y)&\geq \vgf\gd_y -h_y - L_\mu(cg_y)\geq \vgf \gd_y
\EAL \end{equation}
in $B^y_b$ for $b\geq b_0$. This proves \eqref{I.2.1}.
\qed

\noindent\textbf{Completion of proof of Theorem \ref{G-estA}.} \\[3mm]
\underline{The case $N>3$.}\  In $B^y_b$: $|x-y|^{3-N}/\gd(y)<|x-y|^{2-N}/b$. Therefore,  using  \eqref{H-vgf},
\begin{equation}\label{I.2.6}\BAL
 &\vgf_\gV(y)\Gg_y\leq (\vgf_\gV\Gg_y+cg_y)\\
 \leq (&\vgf_\gV(x)+(c/b)\vgf_\gV(y))\Gg_y\leq (c'(b)+c/b) \vgf_\gV(y)\Gg_y \qtxt{in } B^y_b.
\EAL\end{equation}
\medskip

For $a\in \BBR$ and $c$ as in \eqref{I.2.1}
\begin{equation}\label{Fay}
 F_{a,y}:=a\vgf_\gV(y)G^\gV_y-(\vgf_\gV\Gg_y+cg_y).
\end{equation}
Note that, as a distribution, $\vgf_\gV\gd_y= \vgf_\gV(y)\gd_y$. Therefore, by \eqref{I.2.1}, for every
$a\in (0,1)$,
\begin{equation}\label{I.2.7}
-L^\gV(F_{a,y})\leq (a-1)\vgf_\gV(y)\gd_y <0 \qtxt{in $B^y_b\sms\{y\}$}.
\end{equation}

 By Lemma \ref{I.1},  
 for every $a\in(0,1)$ there exists $\gb_y\in(b,\infty)$  \sth
 \begin{equation}\label{I.2.8>}
 F_{a,y}\leq 0 \qtxt{in $B_\gb^y$}, \quad \gb_y <\gb.
\end{equation}
Furthermore, by Lemma \ref{AA'1}, if $a$ is sufficiently
small (depending on $b$ but not on $y$) then
\begin{equation}\label{I.2.9>}
 F_{a,y}\leq 0 \qtxt{on $\prt B^y_b$}.
\end{equation}
Hence for $a$ and $\gb$ as above
$$(F_{a,y})_+=0 \qtxt{on }  \prt B_b^y\cup \prt B_\gb^y$$
and, by  \eqref{I.2.7}, $(F_{a,y})_+$ is $L^\gV$ subharmonic in $B_b^y\sms B_\gb^y$. Consequently
$(F_{a,y})_+=0$ in this domain. As $\gb$ can be chosen arbitrarily large it follows that
$(F_{a,y})_+=0$ in $B_b^y\sms\{y\}$. Therefore, by \eqref{H-vgf}, \eqref{I.2.6}, \eqref{Fay},
\begin{equation}\label{G<Gamma}
 G^\gV_y \leq C_b \Gg_y \qtxt{in}\; B_b^y,
\end{equation}
for $b>b_0$ and $C_b$ independent of $y$. Finally applying Lemma  \ref{b<1} we conclude that \eqref{G<Gamma} holds for every $b>0$.

To obtain the estimate from below, we consider the Green kernel of $L^{\gV}$ in $B_b^y$.
 Clearly
\begin{equation}\label{Gyb}
  G^\gV_{B^y_b}<G^\gV_\Gw \qtxt{in }B^y_b.
\end{equation}

We blow up the ball $B^y_b$ by the transformation:
$\gx= b(x-y)/\gd(y)$ which maps $B^y_b$ to the unit ball $|\gx|<1$. Under this transformation $L^\gV$ becomes,
$$\tl L^W= \Gd_\gx+ W, \qtxt{where}\; |W|\leq \bar a(b-1)^{-2}$$
and $\bar a $ as in \eqref{Vcon1}.
If $\tl G^W$ denotes the Green kernel of $\tl L^W$ in $|\gx|<1$ then,
$$ \tl G^W_0(\gx)= \big(\frac{b}{\gd(y)}\big)^{2-N}G^\gV_{B^y_b}(x,y).$$
It is known that,
$$\rec{c}|\gx|^{2-N}\leq \tl G^W_0(\gx)\leq c |\gx|^{2-N},\quad |\gx|<1/2$$
where $c$ depends only on the bound for $|W|$. Therefore,
$$\rec{C}|x-y|^{2-N}\leq G^\gV_{B^y_b}(x,y) \leq C |x-y|^{2-N},\quad |x-y|<\gd(y)/2b$$
where $C$ depends only on $\bar a $ and $b$. This inequality and \eqref{Gyb} imply
\begin{equation}\label{G>Gamma}
 G^\gV_y \geq C'_b \Gg_y \qtxt{in}\; B_{2b}^y,
\end{equation}
where $C'_b$ is independent of $y$.

\medskip

\noindent\underline{The case $N=3$. }\ If $N=3$, $g_y=\vgf_\gV(y)/\gd(y)$ and $\Gg_y(x)=a_3|x-y|^{-1}$. As before we choose the constant $c$ in \eqref{I.2.1} as follows:
\begin{equation}\label{N=3.a}
 c= C_0/\ell = -C_0 (b-1)^2/\mu,
\end{equation}
where $C_0$ is the constant in \eqref{h_y}.
 Since $c<0$ the completion of the proof requires certain modifications.
\medskip

The constant $c'(b)$ in \eqref{H-vgf} decreases as $b$ increases and $c'(b)\downarrow 1$
as $b\uparrow \infty$. Choose $b_0$ \sth
$$1\leq c'(b)\leq 2 \forevery b>b_0.$$
Put
$$D^y_b:=\{x:|x-y|<\frac{a_3\mu}{4C_0b^2}\gd(y)\}.$$

If   $x\in D^y_b$ then,
$$\frac{C_0b^2}{\mu\gd(y)}<\frac{a_3}{4}|x-y|^{-1}= \rec{4} \Gg_y(x).$$
If, in addition, $ b>b_0$ then,
\begin{equation}\label{N=3.b}\BAL
&\vgf_\gV(x)\Gg_y(x) + cg_y(x)> \vgf_\gV(x)\Gg_y(x) - \frac{C_0b^2}{\mu\gd(y)}\vgf_\gV(y)\\
\geq &\vgf_\gV(y)\big(\rec{c'(b)}\Gg_y(x)- \frac{C_0b^2}{\mu\gd(y)}\big)
\geq \rec{4} \Gg_y(x)\vgf_\gV(y).
\EAL\end{equation}

Let $F_{a,y}$ be as in \eqref{Fay}.  In view of \eqref{N=3.b} we can proceed as before and -- replacing $B^y_b$ by $D^y_b$ -- we obtain,
\begin{equation}\label{G<Gamma'}
 G^\gV_y \leq C_b \Gg_y \qtxt{in}\; D_b^y, \quad b>b_0
\end{equation}
where $C_b$ is independent of $y$. Applying Lemma  \ref{b<1} we conclude that this inequality holds for every $b>0$.

Finally the proof of inequality \eqref{G>Gamma} applies, without modification, to the case $N=3$.
\qed

\section{Proof of Theorem \ref{G-estB}.}

The core of the proof is in the following result.

\begin{lemma}\label{IV.1}
Let $x,y$ be points in $\Gw$ \sth
\begin{equation}\label{IV.1.1}
 |x-y|<r_0/4\kappa, \quad \max (\gd(x),\gd(y))< \frac{|x-y|}{b},\quad b = 16(\kappa+1)^2.
\end{equation}
Then there exists a constant $C'$ depending only on $r_0,\kappa$ and $\bar a$ \sth \eqref{G-est.2} holds.
\end{lemma}

\proof Let $X,Y\in \bdw$ be points \sth
\begin{equation}\label{XY.bdw}
  |x-X|=\gd(x),\quad |y-Y|=\gd(y).
\end{equation}

Let $\gx^X$ be an admissible set of local coordinates at $X$ associated with the cylinder
$T^X(r_0,10\kappa r_0)$ (see \eqref{TP}). Put
$$r=|x-y|,\quad r_1=|\xi^X_1(x-y)|,\quad r'= |(\xi^X)'(x-y)|$$
and similarly $R=|X-Y|$ etc.

The relation
$$\stackrel{\longrightarrow}{xy}\,= \,\stackrel{\longrightarrow}{xX}+ \stackrel{\longrightarrow}{XY}
+\stackrel{\longrightarrow}{Yy}$$
together with \eqref{XY.bdw} and \eqref{IV.1.1} yields
\begin{equation}\label{r-R}\BAL
\max( |r-R|, |r'-R'|, |r_1-R_1|)\leq 2r/b.
\EAL\end{equation}
In particular $R\neq 0$.
 By assumption, (see \eqref{TP}),
 \begin{equation}\label{R1R'}
 R_1=|\xi_1^X(X)- \xi_1^X(Y)|= |f^X((\xi^X)'(X))-f^X((\xi^X)'(Y))|\leq \kappa R'.
 \end{equation}
Hence, by \eqref{r-R}
$$r_1\leq \kappa R' + 2r/b\leq \kappa(r'+ 2r/b) +2r/b= \kappa r' + 2(\kappa+1)r/b.$$
Therefore, with $b$ as in \eqref{IV.1.1},
$$\BAL  r^2&\leq (\kappa r' +2(\kappa+1)r/b)^2+(r')^2\\ &\leq (1+\kappa^2)(r')^2 + 4\kappa(1+\kappa)r^2/b
+ 4(\kappa+1)^2(r/b)^2\\
&\leq (1+\kappa^2)(r')^2 +\frac{r^2}{2}.
\EAL$$
Thus,
\begin{equation}\label{r<gbr'}
 r'\leq r\leq \gb r' \quad \textrm{where} \quad \gb:=\sqrt{2(1+\kappa^2)}.
\end{equation}
It follows that
\begin{equation}\label{x-X}
\gb/b<1/16, \quad \max( |x-X|,|y-Y|)\leq \frac{r}{b}<\frac{r'}{16}.
\end{equation}
Hence, by \eqref{r-R} and \eqref{R1R'}
\begin{equation}\label{r&R}
R' \leq R \leq \sqrt{1+\kappa^2}R', \quad 7r'/8  \leq R' \leq 9r'/8. 
9r_1/8.
\end{equation}
Therefore,
\begin{equation}\label{x_in_TX}\BAL
 x\in T^{X}(3R'/8, 6\kappa R'), \quad y\in T^Y(3R'/8, 6\kappa R')
\EAL\end{equation}
where $T^Y$ is expressed in the coordinates $\xi^Y:=\xi^X- \xi^X(Y)$. Recall,
\begin{equation}\label{Retc}
 |\xi^X(Y)|=|X-Y|=R, \; |(\xi^X)'(Y)|=R',\;  |\xi^X_1(Y)|=R_1\leq \kappa R'.
\end{equation}
Furthermore, as $R'\geq 7r'/8$, we have
\begin{equation}\label{TXTY}
x\nin T^{Y}(3R'/4,6\kappa R'), \quad y\nin T^{X}(3R'/4, 6\kappa R').
\end{equation}
However $T^{X}(3R'/4, 6\kappa R')\cap T^{Y}(3R'/4, 6\kappa R')\cap\Gw \not=\emptyset$. In fact, if
$\Pi_1(X,Y)$ is the half plane whose boundary is the $\xi^X_1$ axis and contains the point $Y$ then,
\begin{equation}\label{=S}\BAL
  \Pi_1(X,Y)\,\cap\, [\xi^{X}_1=4\kappa R']\,\cap\, [|(\xi^{X})'|=R'/2]\\  \sbs T^{X}(3R'/4, 6\kappa
  R')\,\cap \, T^{Y}(3R'/4, 6\kappa R')\,\cap\,\Gw.
\EAL\end{equation}
The intersection on the left hand side consists of a single point  $S$ where
$$\xi^X_1(S)=4\kappa R', \quad |(\xi^{X})'(S)|=R'/2.$$

Applying the BHP theorem in $T^{X}(3R'/4, 6\kappa R')$ when $A_X:=S$ and $A'_X$ is defined by $\xi^X(A'_X):= (5\kappa R', 0)$
yields,
 \begin{equation}\label{GVsim}\BAL
\frac{ G^\gV(x,x_0)}{G^\gV(S,x_0)}&\sim \frac{ G^\gV (x,A'_X)}{G^\gV(S,A'_X)}\\
\frac{ G^\gV(x,x_0)}{G^\gV(S,x_0)}&\sim \frac{ G^\gV(x,y)}{G^\gV(S,y)}
\quad\forevery x\in T^X(3R'/8, 3\kappa R').
\EAL\end{equation}
The second relation is valid because $y\nin T^{X}(3R'/4, 6\kappa R')$.

By Proposition \ref{AA'1}, $G^\gV(S,A'_X)\sim (R')^{2-N}\sim r^{2-N}$. (The relation $r\sim R'$ follows from
\eqref{r<gbr'} and \eqref{r&R}). Moreover, $ G^\gV(x,x_0)\sim \vgf_\gV(x)$. Hence,
\begin{equation}\label{GVsim1}\BAL
G^\gV(x,A'_X)&\sim \frac{\vgf_\gV(x)}{\vgf_\gV(S)}r^{2-N}\\
G^\gV(x,y)&\sim \frac{\vgf_\gV(x)}{\vgf_\gV(S)} G^\gV(S,y)\quad \forevery x\in T^X(3R'/8, 3\kappa R').
\EAL\end{equation}
In these relations the similarity constants depend only on $x_0$, $r_0$, $\kappa$ and $\bar a$.

Next we apply the BHP theorem in $T^{Y}(3R'/4, 6\kappa R')$ when $A'_Y:=S$. Since
$\xi^Y(S)=\xi^X(S)-\xi^X(Y)$ and, by \eqref{R1R'}, $|\xi^X_1(Y)|= R_1\leq \kappa R'$, it follows that
$$3\kappa R'\leq \xi^Y_1(S)\leq 5\kappa R',  \quad |(\xi^{Y})'(S)|=R'/2.$$
We choose $A_Y$ so that $\xi^Y(A_Y)= (2\kappa R',0)$. As in the first relation of \eqref{GVsim1}, we have
\begin{equation}\label{GVsim2}
 G^\gV(y,S)\sim \frac{\vgf_\gV(y)}{\vgf_\gV(A_Y)}r^{2-N} \quad \forevery y\in T^{Y}(3R'/8, \frac{3\kappa}{2}
 R').
\end{equation}
Combining \eqref{GVsim2} and the second relation of \eqref{GVsim1} we obtain,
\begin{equation}\label{GVsim3}
 G^\gV(x,y)\sim \frac{\vgf_\gV(x)}{\vgf_\gV(A_X)}\frac{\vgf_\gV(y)}{\vgf_\gV(A_Y)}r^{2-N}.
\end{equation}
Here we used the symmetry of $G^\gV$ and substituted $S=A_X$. Again the similarity constants depend only on
$x_0$, $r_0$, $\kappa$ and $\bar a$.

Note that $A_X$, $A_Y$ are points lying 'above' $X$ and $Y$ respectively, i.e. on an approximate normal
from the respective boundary point, at a distance proportional to $|x-y|$ which in turn is proportional
to
$|X-Y|$. Applying Lemma \ref{b<b_0} to $\vgf_\gV$ we see that \eqref{GVsim3} remains valid if $A_X, A_Y$ are replaced by any two points in $A(x,y)$. The similarity constant is independent of $r$, but depends on the
proportionality constants mentioned above and therefore on $x_0$, $r_0$, $\kappa$ and $\bar a$.

\qed

\begin{lemma}\label{IV.2}
Let $x,y\in \Gw$ satisfy,
\begin{equation}\label{mixed}
   \min(\gd(x),\gd(y))\leq\frac{|x-y|}{b}\leq \max(\gd(x),\gd(y)), \quad b = 16(\kappa+1)^2.
\end{equation}
Then there exists a constant $C'$ depending only on $r_0,\kappa$ and $\bar a$ \sth \eqref{G-est.2} holds.
\end{lemma}

\proof We assume,
\begin{equation}\label{r/b}
 \gd(x)\leq \frac{|x-y|}{b}< \gd(y).
\end{equation}

Denote by $X$ the point on $\bdw$ \sth $|x-X|=\gd(x)$. Let $\gx^X$ be an admissible set of local
coordinates at $X$ associated with the cylinder $T^X(r_0,10\kappa r_0)$ (see \eqref{TP}). Put
$$r=|x-y|,\quad r_1=|\xi^X_1(x-y)|,\quad r'= |(\xi^X)'(x-y)|.$$
Let $Y\in\bdw$ be the point  \sth $(\xi^X)'(Y)=(\xi^X)'(y)$ and let $R=|X-Y|$ etc.
We consider  the following two cases separately:
\medskip

(a) $r> 2r_1$ \  (b) $r\leq 2r_1$.
\medskip

\ntxt{Case $(a)$} The Lipschitz property of $\bdw$ implies
\begin{equation}\label{base1}
  R_1\leq \kappa R' ,\quad  R'<R<\sqrt{1+\kappa^2}R'.
\end{equation}
Assumption \eqref{r/b} together with the definition of $Y$ imply,
\begin{equation}\label{base2}
  |\xi^X(x)|\leq r/b,\quad |R'-r'|= |(\xi^X)'(X-x)|<r/b.
\end{equation}
As $r-r_1<r'$,  (a) implies that $r'<r<2r'$
and consequently, by \eqref{base2},
\begin{equation}\label{base3}
  r'(1-\frac{2}{b})< R'<r' (1+\frac{2}{b})
\end{equation}

These in turn imply that $|\xi^X(x)|\leq \frac{2}{b-2}R'<R'/32$.
Therefore \eqref{x_in_TX}, \eqref{TXTY} and \eqref{=S} hold and the continuation of the proof is the same
as in the proof of Lemma \ref{IV.1}.
\medskip

\ntxt{Case $(b)$} Let $X$ and $\xi^X$ be as in part (a). By assumption
$r_1=\xi^X_1(y)- \xi^X_1(x)>r/2$ and by construction
$\xi^X_1(x)=\gd(x)<r/b.$
Consequently
\begin{equation}\label{xiy}
  r/2+\gd(x)< \xi^X_1(y)=r_1+\xi^X_1(x)\leq r(1+\rec{b}).
\end{equation}
Moreover, as $(\xi^X)'(x)=0$,
\begin{equation}\label{xi'y}
  (\xi^X)'(y)=r'<(\sqrt{3}/2)r
\end{equation}
We apply Theorem BHP in the standard cylinder $T^X(r,10\kappa r)$. Let $A'=y$ and let $A\in \Gw$ be the
point $\xi^X(A)=(r/4,0)$. Put $v=G^\gV(\cdot,y)$ and $w=G^\gV(\cdot,x_0)$ where $x_0$ is a reference point in
$\Gw$ \sth $\gd(x_0)>r_0$. Then, by BHP,
\begin{equation}\label{vwz}
 \frac{v(z)}{w(z)}\sim \frac{v(A)}{w(A)}  \forevery z\in T^X(r/2, 5\kappa r).
\end{equation}
Recall that $v(A)=G^\gV(A,A')\sim r^{2-N}$ and $w\sim \vgf_\gV$ in $T^X(r,10\kappa r)$. Therefore \eqref{vwz}
implies:
\begin{equation}\label{xyr}
 \rec{C}\frac{\vgf_\gV(x)}{\vgf_\gV(A)}r^{2-N}\leq G^\gV(x,y)\leq C \frac{\vgf_\gV(x)}{\vgf_\gV(A)}r^{2-N}.
\end{equation}
The constant $C$ and all the similarity constants depend only on $x_0$, $\kappa$, $N$ and $\bar a$.

Since $\gd(y)\geq r(1/2-1/b)$ while $\gd(x)\leq r/b$ it follows that $$|x-y| \sim \gd(y)\sim r.$$
Therefore,  by Lemma \ref{b<b_0}, 
$\vgf_\gV(A)\sim \vgf_\gV(x_y)$ for $x_y\in A(x,y)$ so that \eqref{xyr} is equivalent to \eqref{G-est.2}.
\qed

\section{Theorems \ref{K-est} and \ref{3G}}

\subsection{Proof of Theorem \ref{K-est}}

It is well-known that, under the assumptions of the theorem,

\begin{equation}\label{GtoK}
 K^{\gg V}(x,y)=\lim_{z\to y}\frac{G^{\gg V}(x,z)}{G^{\gg V}(x_0,z)}.
\end{equation}
Using the estimate of the Green function \eqref{G-est.2} we obtain (for $z$ near to $y$),

$$\BAL  \frac{G^{\gg V}(x,z)}{G^{\gg V}(x_0,z)}&\sim |x-z|^{2-N}\frac{\vgf_\gV(x)\vgf_\gV(z)}{\vgf_\gV(x_z) \vgf_\gV(z_x)}
\rec{\vgf_\gV(z)}\\ &= |x-z|^{2-N}\frac{\vgf_\gV(x)}{\vgf_\gV(x_z)\vgf_\gV(z_x)}
\EAL $$
where $x_z$ and $z_x$ can be chosen arbitrarily from the set $A(x,z)$ (see \eqref{Axy}). As $z\to
y$ we may replace these points by a point $x_y\in A(x,y)$. This yields \eqref{K-est.1}.

\qed

\subsection{Proof of Theorem \ref{3G}} \ 

In this section $b$ stands for a number that may vary in the interval $[b_0, 2b_0)$, $b_0=\max(b_1,16(1+\kappa)^2)$ and $b_1$ as in \eqref{q-mon}.
   
\textbf{Step 1.} Suppose that each of the pairs $(x,y), (y,z), (x,z)$ satisfies \eqref{d>>r} for some fixed $b>0$. Then, by 
Theorem~\ref{G-estA} inequality \eqref{3G.1} reduces to
\begin{equation}\label{xyz1}
 (|x-y||y-z|)^{2-N}\leq C_b|x-z|^{2-N}\big(|x-y|^{2-N}+|y-z|^{2-N}\big).
\end{equation}

This inequality is easily verified. By the triangle inequality, $$\rec{2}|x-y|\leq\max (|x-z|, |y-z|).$$
 If $|x-y|\leq 2|x-z|$  then
$ (|x-z|/ |x-y|)^{2-N}\leq 2^{N-2}$ and \eqref{xyz1} follows. If $2|x-z|<|x-y|$ then $|x-y|\leq 2|y-z|$ so that
$ (|x-y|/ |y-z|)^{2-N}\leq 2^{N-2}$ and again \eqref{xyz1} follows.
\medskip

\textbf{Step 2.} Assume that the pair $x,z$ satisfies \eqref{d>>r}.
Then, by Theorem~\ref{G-estA},
$$G(x,z)\sim |x-z|^{2-N}.$$
Therefore, by Theorem \ref{G-estB} and inequality \eqref{xyz1}, \eqref{3G.1} reduces to
\begin{equation}\label{3G.2}\BAL
\vgf_\gV(y)^2\vgf_\gV(x)\vgf_\gV(z)\lesssim \vgf_\gV(x_y)^2\vgf_\gV(y_z)^2.
\EAL\end{equation}
where $x_y$ is an arbitrary point in $A(x,y)$ and  $y_z\in A(y,z)$.
We proceed to prove \eqref{3G.2}.

\textbf{Case 2a.} Assume that the pair  $x,y$ satisfies \eqref{d>>r}. Then $x,y\in A_b(x,y)$ and we may choose $x_y=x$ as well as
$x_y=y$. Consequently, 
\begin{equation}\label{d-con1}
\vgf_\gV(x)\vgf_\gV(y) \sim \vgf_\gV(x_y)^2.
\end{equation}
By the same reasoning, if the pair $y,z$ too satisfies \eqref{d>>r} then
\begin{equation}\label{d-con2}
\vgf_\gV(z)\vgf_\gV(y) \sim \vgf_\gV(y_z)^2
\end{equation}
and \eqref{3G.2} holds. 

Now we have to verify \eqref{3G.2} when $y,z$ does not satisfy \eqref{d>>r}, i.e.,
\begin{equation}\label{d-con}
\min (\gd(y),\gd(z))\leq \rec{b}|y-z| 
\end{equation}
for some $b\geq 16(1+\kappa)^2$. We shall show that \eqref{d-con} implies 
\begin{equation}\label{3G.3}\BAL
\vgf_\gV(y)\vgf_\gV(z) \lesssim \vgf_\gV(y_z)^2.
\EAL\end{equation}
We verify this inequality in each of the cases:
\[
\begin{cases}
	(i)\quad  \gd(z)\leq \frac{|y-z|}{b} \leq \gd(y)\\
	(ii)\;\;  \gd(y)\leq \frac{|y-z|}{b} \leq \gd(z)\\
	(iii)\; \max (\gd(y),\gd(z))\leq \rec{b}|y-z|
\end{cases}
\]
If (i) holds then $y\in A_b(y,z)$. Further we choose a point $\gz\in A_b(y,z)$ \sth:
$$ \gd(\gz)=\hat r(y,z)=|y-z|\vee \gd(y)$$ 
and  the pair $z, \gz$ lies  on a $\gl$ pseudo-normal. Therefore we may choose $y_z=y$ as well as  $y_z=\gz$.   By \eqref{q-mon},
$\vgf_\gV(z)\lesssim \vgf_\gV(\gz)$ and \eqref{3G.3}  holds. 

Clearly, the same conclusion holds if (ii) holds.

Finally if (iii) holds we choose $\gz\in A_b(y,z)$ as above and $\eta\in A_b(y,z)$ in the same way except that now the pair $y, \eta$ lies  on a $\gl$ pseudo-normal. By \eqref{q-mon},
$$\vgf_\gV(z)\lesssim \vgf_\gV(\gz), \quad \vgf_\gV(y)\lesssim \vgf_\gV(\eta).$$
Choosing once $y_z=\gz$ and once $y_z=\eta$ we obtain \eqref{3G.3}.

\textbf{Case 2b} Assume that the pair $(y,z)$ satisfies \eqref{d-con} and the pair $x,y$ satisfies a similar inequality:
\[\min (\gd(x),\gd(y))\leq \rec{b}|x-y|.  \]
Then as shown in Case 2a:
\[\vgf_\gV(y)\vgf_\gV(z) \lesssim \vgf_\gV(y_z)^2, \quad \vgf_\gV(x)\vgf_\gV(y) \lesssim \vgf_\gV(x_y)^2\]
\medskip
which implies \eqref{3G.2}.
\medskip

\textbf{Step 3.} \  It remains to consider the case when   the pair  $x,z$ does not satisfy \eqref{d>>r}, i.e.,
\begin{equation}\label{d-con'}
\min (\gd(x),\gd(z))\leq \rec{b}|x-z| 
\end{equation}
for some $b\geq 32  (1+\kappa)^2$. Then, by Theorem \ref{G-estB} and inequality \eqref{xyz1}, \eqref{3G.1} reduces to
\begin{equation}\label{3G.4}\BAL
\vgf_\gV(y)^2\vgf_\gV(x_z)^2\lesssim \vgf_\gV(x_y)^2\vgf_\gV(y_z)^2.
\EAL\end{equation}

\textbf{Case 3a.} \ Assume that,
\begin{equation}\label{max(x,z)}
\max (\gd(x),\gd(z))\leq \rec{b}|x-z|. 
\end{equation}
By the triangle inequality,
$$\max(|x-y|, |y-z|)> \rec{2}|x-z|.$$
Without loss of generality  we assume that the maximum is $|x-y|$ so that 
\begin{equation}\label{xy-dom}
	|x-y|> \rec{2}|x-z|\geq \frac{b}{2}\max (\gd(x),\gd(z)).
\end{equation}

Let $x_y\in A_b(x,y)$ and $x_z\in A_b(x,z)$ be points lying on a  $\gl$ pseudo normal \sth 
$$\gd(x_y) = 2|x-y|, \quad \gd(x_z)= \frac{2}{b}|x-z|.$$ 
(Note that $x$ or $y$  need not be in $A_b  (x,y)$.) 
In view of \eqref{xy-dom} such a choice is possible and by \eqref{q-mon} 
\begin{equation}\label{x_y,x_z}
\vgf_\gV(x_z) \lesssim \vgf_\gV(x_y)    
\end{equation}

Let $y_z\in A(y,z)$ be a point  \sth:\\ 
(i) if $\gd(y)>|y-z|/b$ then 
$y_z=y$. (Note that, as $\gd(z)\leq (1+b)\gd(y)$, $\gd(y)>\hat r(y,z)/(b+1)$).\\ 
(ii) if    $\gd(y)\leq |y-z|/b$ then $\gd(y_z)=|y-z|$ and $y,y_z$ lie on a $\gl$ pseudo normal. 

In either case, using \eqref{q-mon},  we obtain 
\begin{equation}\label{y,y_z}
\vgf_\gV(y) \lesssim \vgf_\gV(y_z)
\end{equation}
Inequality \eqref{3G.4}  follows from \eqref{x_y,x_z} and \eqref{y,y_z}.

\textbf{Case  3b.} \ Assume that:
\begin{equation}\label{min(x,z)}
\gd(x)\leq \rec{b}|x-z|\leq \gd(z).
\end{equation}
Then $\gd(z)\geq \rec{b}\hat r(x,z)$ so that
$z\in A_b(x,z)$. Therefore we may and shall choose $x_z=z$.

Next we choose   
$y_z\in A_b(y,z)$ as follows:
If  $\gd(y)\geq \rec{b}(\gd(z)\vee |y-z|)$ then $\gd(y)\geq \rec{b}\hat r(y,z)$ and we choose  $y_z$ \sth $\gd(y_z)=\gd(y)$ and $z,\,y_z$ lie on a $\gl$ pseudo normal. 

If $\gd(y)\leq \rec{b}(\gd(z) \vee |y-z|)$ then $\gd(y)\leq \rec{b}\hat r(y,z)$. In this case we choose a point $y_z\in A_b(y,z)$ \sth $\gd(y_z)= \hat r(y,z)$ and $z,\, y_z$ lie on a $\gl$ pseudo normal. In either case,
by \eqref{q-mon}, 
$$\vgf_\gV(z)\lesssim \vgf_\gV(y_z).$$

Since, by coice, $z=x_z$ we have  $\vgf_\gV(x_z)\lesssim \vgf_\gV(y_z).$
In order to esablish \eqref{3G.4}, it remains to show that,
$$\vgf_\gV(y)\lesssim \vgf_\gV(x_y)$$
where $x_y$ is a point in $A(x,y)$.
This is proved in the same way as \eqref{y,y_z} replacing $y,\, y_z$ by $y,\, x_y$.


\qed

\noindent
\textbf{Acknowledgement.} The author is grateful to Alano Ancona for several  useful suggestions and to Yehuda Pinchover for many helpful discussions 
on subjects related to this paper.

\end{document}